\documentclass[12pt,letterpaper,reqno]{amsart}

\usepackage{times} 
\usepackage[T1]{fontenc} 
\usepackage{mathrsfs} 
\usepackage{latexsym} 
\usepackage[dvips]{graphics}
\usepackage{epsfig}
\usepackage{amsmath,amsfonts,amsthm,amssymb,amscd} 
\input amssym.def 
\input amssym.tex

\newcommand\be{\begin{equation}} 
\newcommand\ee{\end{equation}}
\newcommand\bea{\begin{eqnarray}} 
\newcommand\eea{\end{eqnarray}} 
\newcommand\bi{\begin{itemize}}
\newcommand\ei{\end{itemize}} 
\newcommand\ben{\begin{enumerate}} 
\newcommand\een{\end{enumerate}}
\newcommand\bc{\begin{center}} 
\newcommand\ec{\end{center}} 
\newcommand\ba{\begin{array}} 
\newcommand\ea{\end{array}}

\theoremstyle{definition} 

\begin{document}

\title[Permutations all of whose patterns of a given length are distinct] 
{Permutations all of whose patterns of a given length are distinct}


\author{Peter Hegarty} \address{Department of Mathematical Sciences, 
Chalmers University Of Technology and University of Gothenburg,
41296 Gothenburg, Sweden} \email{hegarty@chalmers.se}


 \keywords{}

\date{\today}

\begin{abstract} 
For each integer $k \geq 2$, let $F(k)$ denote the largest $n$ for which
there exists a permutation $\sigma \in S_n$ all of whose patterns of length
$k$ are distinct. We prove that $F(k) = k + \lfloor \sqrt{2k-3} \rfloor + \epsilon_k$, where $\epsilon_k \in \{-1,0\}$ for every $k$. Suggestions for further investigations along these lines are discussed. 
\end{abstract}


\maketitle

\setcounter{equation}{0}

\setcounter{equation}{0}

\setcounter{section}{-1}

\section{Notation}

If $f,g : \mathbb{N} \rightarrow \mathbb{R}_{+}$ are two functions, we write
$f(x) \lesssim g(x)$, or alternatively $g(x) \gtrsim f(x)$, to denote that $\limsup_{n \rightarrow \infty} \frac{f(x)}{g(x)} \leq 1$. 
\par As is usual in combinatorics, we will use interval notation for sets of integers. Hence, for real numbers $a \leq b$, the closed interval $[a,b]$ consists of all integers $n$ such that $a \leq n \leq b$, and so on. The interval $[1,n]$ will be denoted simply by $[n]$. Let $S_n$ denote the symmetric group on $n$ letters. We will consider elements of $S_n$ as bijections $\sigma : [n] \rightarrow [n]$, and use the shorthand $\sigma = \sigma_1 \sigma_2 \cdots \sigma_n$ to denote that $\sigma(i) = \sigma_i$. The number $n$ is called the {\em length} of the
permutation. Some further more specialised notation will be introduced below. 

\section{Introduction}

Let $k,n$ be positive integers with $k \leq n$. If $\sigma \in S_n$ and $\pi \in S_k$, then one says that $\sigma$ contains $\pi$ as a {\em pattern} if there exists a $k$-tuple $(a_1,...,a_k)$, with $1 \leq a_1 < a_2 < \cdots < a_k \leq n$ and 
\be
{\hbox{sign}}(\sigma_{a_j} - \sigma_{a_i}) = {\hbox{sign}}(\pi_j - \pi_i), 
\;\;\; {\hbox{for all $1 \leq i < j \leq k$}}.
\ee
Later on, we will need a special notation to distinguish formally between restrictions of $\sigma$ and the patterns they represent, so we may as well introduce this notation now. Let $(a_1,...,a_k)$ and $(a^{\prime}_1,...,a^{\prime}_k)$ be two $k$-tuples as above, and let $\sigma^1$, $\sigma^2$ be the corresponding restrictions of $\sigma$, which are partial functions from $[n]$ to $[n]$. Thus, for example, $\sigma^1$ is the function from $\{a_1,...,a_k\}$ to $[n]$ such that $\sigma^{1}(a_i) = \sigma_{a_i}$. We write, somewhat informally, $\sigma^{1} = \sigma_{a_1} \sigma_{a_2} \cdots \sigma_{a_k}$, and similarly for $\sigma^2$. We shall use the notation $\sigma^1 \subseteq \sigma$ to indicate that $\sigma^1$ is a restriction of $\sigma$. If $\pi \in S_k$ is the pattern represented by $\sigma^1$, we will abuse notaton slightly and also write $\pi \subseteq \sigma$. 
\par The important distinction is the following: we write $\sigma^1 = \sigma^2$ if $a_i = a^{\prime}_i$ for $i = 1,...,k$, whereas we write $\sigma^1 \sim \sigma^2$ if they yield the same pattern, i.e.: if (1.1) holds for both the $a_i$ and the $a^{\prime}_i$ and for the same $\pi \in S_k$. In the latter case one says that $\sigma^1$ and $\sigma^2$ are {\em pattern isomorphic} as restrictions of $\sigma$.  
\\
\\
The study of permutation patterns has developed rapidly over the last 20 years or so: see, for example, the recent book of Kitaev \cite{K} for a comprehensive overview of the literature. Much of the research undertaken is concerned with one or other of two complementary themes :
\\
\\
{\em Pattern avoidance} : Here one is interested in enumerating, as a function of $n$, permutations in $S_n$ which contain no copies of a fixed set of one or more patterns, often of a fixed length $k$. 
\\
{\em Pattern packing} : Here one is interested in constructing permutations which contain as many copies as possible of one or more fixed patterns, or alternatively, which contain as many different patterns as possible. 
\\
\\
There is by now a rather vast literature on the subject of pattern avoidance. Generating-function and other techniques allow for precise enumeration of permutations which avoid specific patterns, and also enable connections to be established to other kinds of permutation statistics. These results are quantitative in nature. Interesting qualitiative results are rarer, but there have been a number of notable achievements: for example, the exponential-growth result of Marcus and Tardos \cite{MT}, and the Kaiser-Klazar theorem \cite{KK} establishing the dichotomy between exponential and polynomial growth. 
\par The literature on pattern packing is smaller but still substantial. For an introduction to the subject of packing copies of a specific pattern, see \cite{AAHHS}. The paper of Miller \cite{M} contains state-of-the-art results on the subject of permutations which contain as many different patterns as possible. It provides the best estimates to date for the following two natural functions :
\par (i) the maximum number pat$(n)$ of possible patterns, of unspecified length, in a permutation of length $n$. Miller proves that 
\be
2^{n} - O(n^2 2^{n - \sqrt{2n}}) \leq {\hbox{pat}}(n) \leq 2^n - \Theta(n 2^{n - \sqrt{2n}}).
\ee
(ii) the minimum length $L(k)$ of a so-called {\em $k$-superpattern}, i.e.: a permutation which contains every $\pi \in S_k$ as a pattern. Since $\left( \begin{array}{c} L(k) \\ k \end{array} \right)
\geq k!$, Stirling's formula gives a trivial lower bound of
\be
L(k) \gtrsim \left( \frac{k}{e} \right)^2.
\ee
Nobody has yet succeeded in improving on this estimate. Miller obtained the best upper bound to date. She exhibited, for every $k$, a $k$-superpattern whose length is at most $k(k+1)/2$. 
\par As Miller remarks in her paper, the problems of estimating the functions in (i) and (ii) above are, loosely speaking, ``dual to one another''. When reading this, it occurred to us to consider the following notion, which seems more directly ``dual'' to the notion of a superpattern :
\\
\\
{\bf Definition 1.1.} Let $k,n$ be natural numbers with $k \leq n$. A permutation $\sigma \in S_n$ is called a {\em $k$-separator} if it contains at most one copy of any $\pi \in S_k$. 
\\
\\
We found no {\em explicit} mention of this concept in the existing literature. The obvious object to study would seem to be the function $F : \mathbb{N} \rightarrow \mathbb{N}$, where $F(k)$ denotes the maximum length of a $k$-separator. The trivial lower bound for $L(k)$ in (1.3) now translates into a trivial upper bound for $F(k)$. For one must have $\left( \begin{array}{c} F(k) \\ k \end{array} \right) \leq k!$ and hence, by Stirling's formula, 
\be
F(k) \lesssim \left( \frac{k}{e} \right)^2.
\ee
However, the property of being a $k$-separator is far more restrictive than this. Indeed, it is almost trivial that 
\be
F(k) < 2k.
\ee
To see this, let $\sigma \in S_n$, $\sigma = \sigma_1 \cdots \sigma_n$, be a $k$-separator. A priori, there are at most $k$ possiblilites for the pattern formed by $\sigma_1 \sigma_2 \cdots \sigma_{k-1} \sigma_T$, as $T$ runs from $k$ up to $n$. Hence, if $n \geq 2k$, at least two of these patterns must coincide. 
\par The main result of our note is the following :
\\
\\
{\bf Theorem 1.2.} {\em 
For each $k \geq 2$ one has
\be
F(k) = k + \lfloor \sqrt{2k-3} \rfloor + \epsilon_k,
\ee
where $\epsilon_k \in \{-1,0\}$.}
\\
\\
The proof of this result, which follows in Section 2, has much in common with the methods of \cite{M}. To obtain a lower bound for $F(k)$, we employ the same ``tilted checkerboard'' permutations appearing in \cite{M}. For the upper bound, we further extend the idea employed in \cite{M}, and attributed originally to Coleman \cite{C}, that to avoid repeating patterns in a permutation $\sigma \in S_n$, the so-called {\em taxicab distance} between elements $i,j \in [n]$ should be large. What seems to be new in our proof is a sort of optimisation argument which allows for a very precise estimate for $F(k)$ (there are hints of this argument in Section 6 of \cite{M}, but our approach seems to be different).    
\par Our paper closes with a short discussion section (Section 3), which includes some suggestions for extending the ideas introduced here.  
   
\setcounter{equation}{0}

\section{Proof of Theorem 1.2}

\subsection{Proof of Lower Bound}

In this subsection we prove that, for every $k \geq 2$, 
\be
F(k) \geq k + \lfloor \sqrt{2k-3} \rfloor - 1.
\ee
Note that $2k-3$ is a perfect square if and only if $k = 2m^2 - 2m+2$ for some 
$m \geq 1$. To prove (2.1) it thus suffices to prove, for every integer $m \geq 1$, the following two statements :
\\
\\
(i) If $2m^2 - 2m + 2 \leq k \leq 2m^2 + 1$, then there exists $\sigma \in S_{k+(2m-2)}$ which is a $k$-separator,
\\
(ii) If $2m^2 + 2 \leq k \leq 2m^2 + 2m + 1$, then there exists $\sigma \in S_{k+(2m-1)}$ which is a $k$-separator. 
\\
\\
For positive integers $r,s$, we employ the definitions of the $r \times s$
{\em titled rectangle} and the $r \times s$ {\em titled checkerboard} as given in \cite{M}{\footnote{Miller already has some nice pictures as visual aids, so we do not reproduce these here.}}. Let $\sigma^{r,s}$ denote the corresponding $r \times s$ checkerboard permutation, considered as a permutation of length $\lceil \frac{rs}{2} \rceil$. Thus, $r$ is the number of columns and $s$ the number of rows in this permutation. If $a \in [r]$ and $b \in [s]$ are such that $a \equiv b \; ({\hbox{mod $2$}})$, then $\sigma^{r,s}_{a,b}$ denotes the element of the $r \times s$ checkerboard which lies in its $a$:th column and $b$:th row. Here, the columns are read from left to right and the rows from bottom to top.
\\
\\
{\sc Case 1: $2m^2 -2m + 2 \leq k \leq 2m^2 + 1$.}
\\
\\
We consider three subcases :
\\
\par (i.a) $k = 2m^2 - 2m+2$.
\par (i.b) $k = 2m^2 - 2m+2 + i$, for some $1 \leq i \leq m$.
\par (i.c) $k = 2m^2 - m + 2 + j$, for some $1 \leq j < m$.
\\
\\
{\em Subcase (i.a):} We need to exhibit $\sigma \in S_{2m^2}$ which is a $k$-separator. We take $\sigma = \sigma^{2m,2m}$. This is indeed a permutation of length $\frac{1}{2} (2m)^2 = 2m^2$. By Proposition 4.4 of \cite{M}, a pattern $\pi \subseteq \sigma$ which is not uniquely represented must truncate or avoid at least two of the rows and/or columns of $\sigma$. But, clearly, any such $\pi$ must omit at least $2m-1$ elements of the checkerboard.
\\
\\
{\em Subcase (i.b):} We need to exhibit $\sigma \in S_{2m^2 + i}$ which is a 
$k$-separator. We take $\sigma$ to be the prefix of $\sigma^{2m+1,2m}$ whose 
complement consists of its last $m+1-i$ elements, i.e.: we omit the elements 
$\sigma^{2m+1,2m}_{2m+1,b}$ for $b \geq 2i+1$. We need to show that any $\pi \subseteq \sigma$ which omits $2m-2$ elements of $\sigma$ is represented uniquely. We already know this is true when $i = 0$ from subcase (i.a) above. Now suppose $i > 0$ and suppose $\pi^1, \pi^2$ are two restrictions of $\sigma$ such that $\pi^1 \sim \pi^2$ and each omits $2m-2$ elements of $\sigma$. If $\pi^1 \neq \pi^2$ then, reading both from left-to-right, there must be a first position where they differ. Say this is in position $\xi \in [k]$ and let $\pi^{1}_{\xi} = \sigma_{a,b}$, $\pi^{2}_{\xi} = \sigma_{c,d}$. Since both $\pi^1$ and $\pi^2$ omit only $2m-2$ elements of the $2m \times 2m$ checkerboard formed by all but the last column of $\sigma$, at least one of $a$ and $c$ must equal $2m+1$. 
\par First suppose $a = c = 2m+1$ and, WLOG, that $b < d$. Then 
\be
\sigma_{2m+1,b} <  \sigma_{e,f} < \sigma_{2m+1,d}, \;\;\; {\hbox{whenever $e < 
2m+1$, $b \leq f < d$ and $e \equiv f \; ({\hbox{mod $2$}})$}}.
\ee
Since $d \geq b + 2$, there are at least $2m$ such elements $\sigma_{e,f}$. 
Since $\pi^1$ and $\pi^2$ coincide to the left of position $\xi$, none of these $2m$ elements of $\sigma$ can lie in either $\pi^j$. But this contradicts the assumption that the $\pi^j$ omit only $2m-2$ elements of $\sigma$. 
\par Now suppose, WLOG, that $2m+1 = c > a$. Then all elements $\sigma_{e,f}$ must be missing from $\pi^2$, where either (i) $e=a,\; f> b$ (ii) $a < e < 2m+1$ (iii) $e = 2m+1, f < d$ (iv) $e < 2m+1,\; d \leq f < b$ or (v) $e < a, \; b \leq f < d$ (note that either (iv) or (v) is unsatisfiable). It is easy to see that the total number of pairs $(e,f)$ of the same parity satisfying at least one of these conditions must then be at least $2m$, again a contradiction. 
\\
\\
{\em Subcase (i.c):} We need to exhibit $\sigma \in S_{2m^2 + 2m + j}$ which is a $k$-separator. We take $\sigma$ to be the restriction of $\sigma^{2m+1,2m+1}$ whose complement consists of the leftmost $m+1-j$ elements in its top row, i.e.: we omit the elements
$\sigma^{2m+1,2m+1}_{2\xi-1,2m+1}$, where $1 \leq \xi \leq m+1-j$. We need to show that any $\pi \subseteq \sigma$ which omits at $2m-2$ elements is represented uniquely. We already know this is true when $j = 0$ from subcase (i.b) above. For $j > 0$, the argument is essentially the same as in subcase (i.b), for we can rotate $\sigma$ by 90 degrees clockwise and apply a symmetry argument. We shall flesh out the details a little so as to leave no room for doubt. Suppose $\pi^1, \pi^2$ are two restrictions of $\sigma$ such that $\pi^1 \sim \pi^2$ and each omits $2m-2$ elements of $\sigma$. Firstly, if $\pi^1$ and $\pi^2$ coincide along the top row of $\sigma$, then it is easy to see that the restrictions of both to the remaining rows must also be pattern isomorphic. Then we can apply subcase (i.b) directly. So we may suppose that, reading from left-to-right, there is a first position along the top row of $\sigma$, say $\sigma_{a,2m+1}$, such that, WLOG, $\sigma_{a,2m+1} \in \pi^1 \backslash \pi^2$. Suppose $\sigma_{a,2m+1}$ is the $t$:th largest element in $\pi^1$. If the $t$:th largest element in $\pi^2$ appears in the top row of $\sigma$, then it must be in the $(a+2)$:nd column or later. Since $\pi^1 \sim \pi^2$, the same number of elements appear in both to the left of the $t$:th largest element. It follows unavoidably that $\pi^2$ omits at least $2m$ elements of $\sigma$, a contradiction. A similar argument can be applied if the $t$:th largest element of $\pi^2$ doesn't appear in the top row. We will unavoidably be led to the contradiction that $\pi^2$ omits at least $2m$ elements of $\sigma$. 
\\
\\
{\sc Case 2: $2m^2 + 2 \leq k \leq 2m^2 + 2m + 1$.}
\\
\\
We consider three subcases :
\\
\par (ii.a) $k = 2m^2 + 2$.
\par (ii.b) $k = 2m^2 + 2 + i$, for some $1 \leq i \leq m$.
\par (ii.c) $k = 2m^2 + m + 2 + j$, for some $1 \leq j < m$.
\\
\\
{\em Subcase (ii.a):} We take $\sigma = \sigma^{2m+1,2m+1}$, which has length
$\lceil \frac{(2m+1)^2}{2} \rceil = 2m^2 + 2m + 1$. We need to show that any restriction which omits at most $2m-1$ elements of $\sigma$ uniquely represents its pattern. Let $\pi^1$ and $\pi^2$ be two such restrictions and suppose $\pi^1 \sim \pi^2$, but $\pi^1 \neq \pi^2$. We can argue as in subcase (i.c) above that, unless $\pi^1$ and $\pi^2$ both contain the entire top row of $\sigma$, at least one of them necessarily omits at least $2m$ elements of $\sigma$, a contradiction. Hence, both contain the entire top row. By a symmetry argument, both must also contain the entire first column. Hence, neither can truncate a row or column of $\sigma$. By Proposition 4.4 of \cite{M}, the only way left for $\pi^1$ and $\pi^2$ to represent the same pattern is for each to omit an even numbered row and column, and nothing else. But because each of $\pi^1$ and $\pi^2$ contains the entire top row, it is then easy to see that we cannot get the same pattern unless they omit exactly the same row and column, and hence are after all equal.
\\
\\
We can deal with subcases (ii,b) and (ii.c) in exactly the same way as we did with (i.b) and (i.c) respectively, except that now we start with $\sigma^{2m+1,2m+1}$ and add in the elements of the last column and top row of $\sigma^{2m+2,2m+2}$ one-by-one, first going from bottom to top in the last column and then from right to left in the top row. As in Case 1, this shows that the value of $F(k) - k$ is non-decreasing for $k$ in the interval covered by Case 2. 
     
\subsection{Proof of Upper Bound}

In this subsection we prove that, for every $k \geq 2$,
\be
F(k) \leq k + \lfloor \sqrt{2k-3} \rfloor.
\ee
The basic idea is that a $k$-separator must take numbers which
are close together and permute them so they are placed far apart. We now make this precise. We need a couple of definitions:
\\
\\
{\bf Definition 2.1.} Let $\sigma \in S_n$ and let $i,j \in [n]$. The {\em distance} between $i$ and $j$ in $\sigma$, denoted $d_{\sigma}(i,j)$, is defined as
\be
d_{\sigma}(i,j) = |\sigma^{-1}(i) - \sigma^{-1}(j)|.
\ee
In other words, $d_{\sigma}(i,j)$ is the number of spaces between $i$ and $j$
in the representation $\sigma = \sigma_1 \sigma_2 \cdots \sigma_n$. For example, if $n = 9$ and $\sigma = 341679825$, then $d_{\sigma}(4,2) = 6$ and $d_{\sigma}(1,9) = 3$.   
\\
\\
{\bf Definition 2.2.} Let $\sigma \in S_n$ and let $i,j \in [n]$. We define the natural number $t^{\sigma}_{i,j}$ to be the number of integers $l$ such that 
\be
{\hbox{sign}}(l-i) = {\hbox{sign}}(j-l)
\;\;\; {\hbox{and}} \;\;\; {\hbox{sign}}(\sigma^{-1}(l) - \sigma^{-1}(i)) = {\hbox{sign}}(\sigma^{-1}(j) - \sigma^{-1}(l)).
\ee
In other words, $t^{\sigma}_{i,j}$ is the number of integers lying strictly between $i$ and $j$ which also appear between $i$ and $j$ in the representation $\sigma = \sigma_1 \sigma_2 \cdots \sigma_n$. For example, if $n = 9$ and $\sigma = 341679825$, then $t^{\sigma}_{3,8} = 3$, since each of the numbers $4,6,7$ appear between $3$ and $8$ in $\sigma$, whereas $t^{\sigma}_{1,5} = 1$, since only $2$ appears between $1$ and $5$ in $\sigma$. Note that, a priori, for any $i,j$ and $\sigma$ one has 
\be
0 \leq t^{\sigma}_{i,j} < |i-j|.
\ee
We now require three lemmas :
\\
\\
{\bf Lemma 2.3.} {\em Let $\sigma \in S_n$ be a $k$-separator. Then for any $i \neq j \in [n]$ one has 
\be
d_{\sigma}(i,j) \geq (n-k+2) + t^{\sigma}_{i,j} - |i-j|.
\ee}

{\em Proof.} Without loss of generality, $i < j$ and $i$ appears to the left of $j$ in the
standard representation of $\sigma$. We consider 
\be
\sigma = \sigma_1 \cdots \sigma_{r-1} \; {\hbox{{\bf i}}} \; \sigma_{r+1} \cdots \sigma_{s-1} \; {\hbox{{\bf j}}} \; \sigma_{s+1} \cdots \sigma_n.
\ee
By definition, $t^{\sigma}_{i,j}$ of the numbers in the interval $(i,j)$ appear among $\sigma_{r+1},...,\sigma_{s-1}$. Hence, the remaining $j - i - 1 - t^{\sigma}_{i,j}$ such numbers appear either to the left of {\bf i} or to the right of {\bf j}
in (2.8). A total of $n - d_{\sigma}(i,j) - 1$ numbers appear either to the 
left of {\bf i} or to the right of {\bf j}. Hence, exactly 
\be
(n - d_{\sigma}(i,j) - 1) - (j-i-1-t^{\sigma}_{i,j}) = n + t^{\sigma}_{i,j} - (j-i) - d_{\sigma}(i,j)
\ee
of these numbers are not in the closed interval $[i,j]$. If (2.7) failed to 
hold, it would mean that the right-hand side of (2.9) was greater than or equal to $k-1$. In other words, it would mean that at least $k-1$ of the numbers appearing either to the left of {\bf i} or to the right of {\bf j} in $\sigma$ were not in the interval $[i,j]$. If so, pick any $k-1$ such numbers reading from left to right, say $\sigma_{i_1}, \sigma_{i_2},...,\sigma_{i_p}, \sigma_{j_1},...,\sigma_{j_q}$, where $p+q = k-1$ and $i_1 < i_2 < \cdots < i_p < r < s < j_1 < \cdots < j_q$. Then the two subsequences 
\be
\sigma_{i_1} \cdots \sigma_{i_p} \; {\hbox{{\bf i}}} \; \sigma_{j_1} \cdots \sigma_{j_q}
\;\;\;\; {\hbox{and}} \;\;\;\;
\sigma_{i_1} \cdots \sigma_{i_p} \; {\hbox{{\bf j}}} \; \sigma_{j_1} \cdots \sigma_{j_q} 
\ee
yield two copies of the same length-$k$ pattern in $\sigma$, contradicting the fact that $\sigma$ is a $k$-separator. This completes the proof of the lemma.
\\
\\
{\bf Remark 2.4.} Let $\sigma \in S_n$ and $i,j \in [n]$. The {\em taxicab distance} between $i$ and $j$ with respect to $\sigma$, which we denote $d^{{\hbox{TC}}}_{\sigma}(i,j)$, is defined in Section 6 of \cite{M} as 
\be
d^{{\hbox{TC}}}_{\sigma}(i,j) = |i-j| + |\sigma_i - \sigma_j|.
\ee
Also, in Millers' notation, 
\be
t^{\sigma}_{i,j} = n - |S_{\sigma}(i,j)|.
\ee
Hence, in Miller's notation, (2.7) becomes 
\be
d^{{\hbox{TC}}}_{\sigma^{-1}}(i,j) \geq 2n+2 - k - |S_{\sigma}(i,j)|.
\ee 
{\bf Lemma 2.5.} {\em Let $m \in \mathbb{N}$ and let $a_1 a_2 \cdots a_m$ be any permutation of the integers in $[m]$. Then 
\be
\sum_{i=1}^{m-1} |a_i - a_{i+1}| \leq \frac{(m-1)(m+1)}{2}.
\ee}

{\em Proof.} Let 
\be
r := \# \{i \in [m-1] : a_{i+1} < a_i\}, \;\;\; s := m-1-r = \# \{i \in [m-1] :
a_{i+1} > a_i \}.
\ee
Let $i_1,...,i_r$ be the indices such that $a_{i_j + 1} < a_{i_j}$, for $j = 1,..,r$. Since the ``take-off points'' $a_{i_j}$ are distinct, for $j = 1,..,r$, and also the ``landing points'' $a_{i_j + 1}$ are distinct, for $j = 1,...,r$, it follows that
\be
\sum_{j=1}^{r} |a_{i_j} - a_{i_j + 1}| \leq r(m-1) - 2 \cdot \sum_{j=1}^{r} (j-1) = r(m-1) - r(r-1).
\ee
Similarly, if $i^{\prime}_1,...,i^{\prime}_s$ are the indices such that $a_{i^{\prime}_j + 1} > a_{i^{\prime}_j}$, for $j = 1,...,s$, then 
\be
\sum_{j=1}^{s} |a_{i^{\prime}_j} - a_{i^{\prime}_j + 1}| \leq s(m-1) - s(s-1).
\ee 
From (2.16) and (2.17) it follows that 
\be
\sum_{i=1}^{m-1} |a_i - a_{i+1}| \leq m(m-1) - (r^2 + s^2).
\ee
Since $r+s = m-1$, the quantity $r^2 + s^2$ is minimised when $m$ is odd and $r = s = \frac{m-1}{2}$. This proves the lemma. 
\\
\\
{\bf Remark 2.6.} A more careful analysis yields 
\be
\sum_{i=1}^{m-1} |a_i - a_{i+1}| \leq \frac{(m-1)(m+1)}{2} - 1,
\ee
which is best-possible. We have no use for this slight improvement in what follows, however.
\\
\\
{\bf Lemma 2.7.} {\em For any $x \in \mathbb{N}$ one has 
\be
\lfloor \sqrt{x} \rfloor + \frac{x}{\lfloor \sqrt{x} \rfloor} \leq 2\left(\lfloor \sqrt{x} \rfloor + 1 \right),
\ee
with equality if and only if $x = n^2 - 1$ for some $n \in \mathbb{N}$.}
\\
\\
{\em Proof.} This is a simple exercise.
\\
\\
We are now ready to prove (2.3). Let $\sigma \in S_n$ be a $k$-separator. Since the numbers $t^{\sigma}_{i,j}$ in (2.5) are, at the very least, non-negative, we have by Lemma 2.3, for any $i \neq j \in [n]$, that
\be
d_{\sigma}(i,j) \geq (n-k+2) - |i-j|.
\ee
Let $m \in [n]$. The value of $m$ will be optimised in due course. Reading from left to right in the representation $\sigma = \sigma_1 \cdots \sigma_n$, the numbers from $1$ to $m$ will appear in some order, say as 
$\sigma_{i_1},...,\sigma_{i_m}$, where $i_1 < \cdots < i_m$. By (2.21), we have 
that 
\be
\sum_{l=1}^{m-1} d_{\sigma} (\sigma_{i_l}, \sigma_{i_{l+1}}) \geq (m-1)(n-k+2) - 
\sum_{l=1}^{m-1} |\sigma_{i_l} - \sigma_{i_{l+1}}|.
\ee
The left-hand side of (2.22) cannot exceed $n-1$, since we are reading along $\sigma$ from left to right. By Lemma 2.5, the sum on the right-hand side cannot exceed $\frac{(m-1)(m+1)}{2}$. Hence,
\be
n-1 \geq (m-1)(n-k+2) - \frac{(m-1)(m+1)}{2},
\ee
which we can rewrite as 
\be
n \leq (k-1) + \frac{1}{2} f_k (m),
\ee
where 
\be
f_k (m) = m + \frac{2k-3}{m-2}.
\ee
Inequality (2.24) must hold for any choice of $m \in [n]$, so we choose $m$ to make $f_k (m)$ as small as possible. As a function of a real variable, $f_k (m)$ has a local minimum at $m = \sqrt{2k-3} + 2$. Since, in our case, $m$ must be an integer, we take 
\be
m = \left\{ \begin{array}{lr} \lceil \sqrt{2k-3} \rceil + 2, &  {\hbox{if $k = 2u^2 + 1$ for some $u \in \mathbb{N}$}}, \\
\lfloor \sqrt{2k-3} \rfloor + 2, & {\hbox{otherwise}}. \end{array} \right.
\ee
With the help of Lemma 2.7 one easily verifies that, for this choice of $m$ one always has  
\be
f_{k}(m) < 2 \left( \lfloor \sqrt{2k-3} \rfloor + 2 \right)
\ee
which, together with (2.24), yields (2.3).
      
\setcounter{equation}{0}

\section{Discussion}

Since $F(2) = 2$ and $F(3) = 4$, one has $\epsilon_2 = -1$ whereas $\epsilon_3 = 0$. We do not know any nice way to determine the value of $\epsilon_k$ in general. Indeed, we do not even know whether both $-1$ and $0$ occur infinitely often. It is natural to guess that any $k$-separator of length $F(k)$ must look very similar to those considered in subsection 2.1, and hence basically look like a tilted checkerboard, minus some elements around the edges. However, we do not know any elegant approach to tackling these refinements of the extremal problem which has been considered in this article.  
\par There are other directions in which one might choose to extend the ideas presented here, so let us just make one of the more obvious suggestions. For each pair $(k,l)$ of positive integers, let $P(k,l)$ denote the maximum number of distinct patterns $\pi \in S_k$ which can appear in a permutation $\sigma \in S_{k+l}$ and set 
\be
Q(k,l) := \frac{P(k,l)}{\left( \begin{array}{c} k+l \\ k \end{array} \right)}.
\ee
Obviously, for any fixed $k$, $Q(k,l) \rightarrow 0$ as $l \rightarrow \infty$. 
The basic question then is, how quickly does $Q(k,l)$ go to zero ? Theorem 1.2 says that $Q(k,l) = 1$ at least for $l \leq k + \lfloor \sqrt{2k-3} \rfloor - 1$. 




\vspace*{1cm}

\begin{thebibliography}{AAHHS} 

\bibitem[AAHHS]{AAHHS} M.H. Albert, M.D. Atkinson, C.C. Handley, D.A. Holton and W. Stromquist, \emph{On packing densities of permutations}, Electron. J. Combin. \textbf{9} (2002), No.1, Research Paper 5, 20pp. 

\bibitem[C]{C} M. Coleman, \emph{An answer to a question by Wilf on packing distinct patterns in a permutation}, Electron. J. Combin. \textbf{11} (2004), No.1, Note 8, 4pp. 

\bibitem[K]{K} S. Kitaev, \emph{Patterns in permutations and words}, Springer (2011). 

\bibitem[KK]{KK} T. Kaiser and M. Klazar, \emph{On growth rates of closed permutation classes}, Electron. J. Combin. \textbf{9} (2002/03), No.2, Research Paper 10, 20pp.

\bibitem[M]{M} A. Miller, \emph{Asymptotic bounds for permutations containing many different patterns}, J. Combin. Theory Ser. A \textbf{116} (2009), No.1, 92--108.

\bibitem[MT]{MT} A. Marcus and G. Tardos, \emph{Excluded permutation matrices and the Stanley-Wilf conjecture}, J. Combin. Theory Ser. A \textbf{107} (2004), No.1, 153--160.

\end{thebibliography}
\end{document}